\documentclass[a4paper,reqno]{amsart}
\usepackage{hyperref}

\usepackage{amssymb,amsxtra}
\usepackage{graphicx} \usepackage{amsmath} \usepackage{amssymb}
\usepackage[svgnames]{xcolor}
\usepackage[color,matrix,arrow]{xy}
\input xy \xyoption{all}

\newtheorem{lem}{Lemma}[section]
\newtheorem{thm}[lem]{Theorem}
\newtheorem{prop}[lem]{Proposition}
\newtheorem{coro}[lem]{Corollary}
\newtheorem{defn}[lem]{Definition}
\newtheorem{example}[lem]{Example}
\newtheorem{rmk}[lem]{Remark}

\newcommand{\sk}[1]{^{(#1)}}
\newcommand{\su}[2]{^{(#1)}_{#2}}

\newcommand{\wto}{\widetilde{\omega}}
\newcommand{\wt}[1]{\widetilde{#1}}

\newcommand{\Hom}{\mathrm{Hom}}

\newcommand{\uhom}[1]{\underline{\mathrm{Hom}}}
\newcommand{\Ext}{\mathrm{Ext}}
\newcommand{\Ker}{\mathrm{Ker}\,}

\newcommand{\add}{\mathrm{add}\,}

\newcommand{\mmod}{\mathrm{mod}\,}

\newcommand{\dmn}{\mathrm{dim}\,}
\newcommand{\hls}{\mathrm{det}\,}
\newcommand{\arr}[2]{\begin{array}{#1}#2\end{array}}
\newcommand{\mat}[2]{\left(\begin{array}{#1}#2\end{array}\right)}

\begin{document}
\title{$n$-complete algebras and McKay quivers}
\author[Zhang ET AL.]{Tongliang Zhang}
\address{College of Mathematics and Computer Science, Key Laboratory of HPCSIP (Ministry of Education of China),\\ Hunan Normal University, Changsha 410081, CHINA }
\author[]{Deren Luo}
\address{College of Mathematics and Computer Science, Key Laboratory of HPCSIP (Ministry of Education of China),\\ Hunan Normal University, Changsha 410081, CHINA }
\author[]{Lijing Zheng$^{\ast}$}
\address{School of Mathematics and Physics \\ University of South China \\ Hengyang, 421001, Hunan, P. R. China}

\thanks{$\ast$ E-mail: zhenglijing817@163.com\\
This work is partly supported by  Natural Science Foundation of China \#11271119 and by Hunan Provincial Innovation Foundation For Postgraduate(CX2013B216 and CX2014B189).\\}

\maketitle
\begin{abstract}
Let $\Gamma^{n}$ be the cone of an $(n-1)$-complete algebra over an algebraically closed field $k$. In this paper, we prove that if the bound quiver $(Q_{n},\rho_{n})$ of $\Gamma^{n}$ is a truncation  from the bound McKay  quiver $(Q_{G},\rho_{G})$ of a finite subgroup $G$ of $GL(n,k)$,  the  bound quiver $(Q_{n+1}, \rho_{n+1})$ of $\Gamma^{n+1}$, the cone of $\Gamma^{n}$,  is a truncation from the bound McKay quiver $(Q_{\widetilde{G}},\rho_{\widetilde{G}})$ of $\widetilde{G}$, where   $\widetilde{G}\cong G\times \mathbb{Z}_{m}$ for some $m\in \mathbb{N}$.
\end{abstract}

Recently, Iyama introduced $n$-cluster tilting subcategories and developed higher Auslander-Reiten theory (\cite{i}).
In \cite{i1}, he introduced and characterized a class of higher representation algebras,  $n$-complete algebras.
Such algebras are preserved under cone constructions.
He also proved that $n$-Auslander absolutely $n$-complete algebra are constructed by iterative cone construction starting from some path algebra of quiver of type $A_{r}$ with linear orientation.
Guo proved that such algebras can be obtained from a truncation of the McKay quivers of some abelian groups (\cite{g-McKay}).

In this paper, we generalize the results of Guo and prove the following result:
Let $\Gamma^{n}$ be the cone of an $(n-1)$-complete algebra,  if the bound quiver $(Q_{n},\rho_{n})$ of $\Gamma^{n}$ is a truncation  from the bound McKay  quiver $(Q_{G},\rho_{G})$ of a finite subgroup $G$ of $GL(n,k)$, then there exists a positve integer $m$ such that the  bound quiver $(Q_{n+1},  \rho_{n+1})$ of $\Gamma^{n+1}$  is a truncation from the bound McKay quiver $(Q_{\widetilde{G}},\rho_{\widetilde{G}})$ of a finite subgroup $\widetilde{G}\cong G\times \mathbb{Z}_{m}$ in $GL(n+1,k)$.

The paper is organized as follows. In Section 1, we shall recall some basic definitions and facts needed for $n$-complete algebras, McKay quivers and trivial extensions of graded self-injective algebras. Then we describe the bound McKay quivers using twisted trivial extensions in Section 2 and Section 3. Our main theorem will be stated and proved in Section 4.

\section{Preliminaries}
Throughout this paper, $k$ is an algebraically closed field of characteristic $0$.
Let $\Lambda$ be an algebra over $k$.
Denote by $\mmod\Lambda$ the category of finitely generated left $\Lambda$-modules, and  for $M\in\mmod\Lambda$, denote by $\add M$ the full subcategory of $\mmod\Lambda$ consisting of direct summands of finite direct sums of copies of $M$, and denote by $D$ the standard $k$-duality $\Hom_k(-,k)$.

\subsection{$n$-complete algebra}

In \cite{i}, Iyama introduced and studied $n$-complete algebra.
Let $\Lambda$ be a finite dimensional $k$-algebra.
Let $\tau_{n}=D{\rm Tr}\Omega^{n-1}$ be the $n$-Auslander-Reiten translation, the subcategory  $\mathcal{M} = \mathcal{M}_{n}( D\Lambda ) = \add \{ \tau^{i}_{n} (D\Lambda) | i\geq 0 \}$  of $\mmod\Lambda$ is called the {\em $\tau_{n}$-closure} of $D\Lambda$.
Let
\begin{eqnarray*}
\mathcal{I}(\mathcal{M})&=&\add  D\Lambda,\\
\mathcal{P}(\mathcal{M})&=&\{X\in\mathcal{M}~|~\tau_{n}X=0\},\\
 \mathcal{M}_{P}&=&\{X\in\mathcal{M}~|~X  \text{ has no non-zero summands
in}\ \mathcal{P}(\mathcal{M})\},\\
 \mathcal{M}_{I}&=&\{X\in\mathcal{M}~|~X  \text{ has no non-zero summands
in}\ \mathcal{I}(\mathcal{M})\}.
\end{eqnarray*}

$\Lambda$ is called {\em $n$-complete} if
the global dimension of $\Lambda$, gl.dim$\Lambda \leq n$, and the following conditions $(1)-(3)$ are
satisfied.

(1). There exists a tilting $\Lambda$-module $T$ satisfying $\mathcal{P}(\mathcal{M})=\add T$,

(2). $\mathcal{M}$ is an $n$-cluster titing subcategory of $T^{\bot}=\{X\in\mathcal{M}|\Ext^{i}_{\Lambda}(T,X)=0(0<i)\}$,

(3). $\Ext^{i}_{\Lambda}(\mathcal{M}_{P},\Lambda)=0$ for any $0<i<n$.

In this case, $ \Gamma: = \mathrm{End}_{\Lambda}( \bigoplus \limits_{t\geq0} \tau_{n}^{t} ( D\Lambda) ) $  is called the {\em cone} of $\Lambda$.

Let $\Lambda$ be an $n$-complete algebra and $\mathcal{M} = \mathcal{M}_{n} (D\Lambda)$ the $\tau_{n}$-closure of
$D\Lambda$, we denote by $J_{\mathcal{M}}$ the Jacobson radical of  $\mathcal{M}$(\cite{elements}).
Define quiver $(Q_{0} = Q_{0}(\mathcal{M}), Q_{1}= Q_{1} (\mathcal{M}))$ and a map $\tau_{n}$ as follows:

(1).  $Q_{0}$ (respectively,$Q_{P}$ and $Q_{I}$) is the set of indecomposable objects in $\mathcal{M}$ (respectively, $\mathcal{M}_{P}$ and $\mathcal{M}_{I}$).

(2). For $X,Y$$\in$$Q_{0}$, put $d_{XY} =\dmn( J_{\mathcal{M}} (X,Y) /J^{2}_{\mathcal{M}} (X,Y) )$ and draw $d_{XY}$ arrows from $X$ to $Y$.

(3). A map $\tau_{n}: Q_{P} \rightarrow Q_{I}$ is given by   $\tau_{n}: \mathcal{M}_{P} \rightarrow \mathcal{M}_{I}$.

$(Q_0, Q_1, \tau_n)$ is a weak translation quiver in the sense of \cite{i} and is  called the {\em Auslander-Reiten quiver of $\mathcal{M}$}.

Let $\Gamma^{n}$ be  the cone of an $(n-1)$-complete algebra $\Gamma^{n-1}$, and let $(Q_{n,0}, Q_{n,1}, \tau_{n})$ be the Auslander-Reiten quiver of $\mathcal{M}_{n-1}(D\Gamma^{n-1})$.
By Lemma 6.4 of \cite{i}, $\Gamma^{n}$ is given by the bound quiver $(Q_n, \rho_n)$, that is $\Gamma^{n} \cong kQ/(\rho_{n})$ for  some relation $\rho_n$.
We say that $\Gamma^{n}$ is given by the {\em bound quiver} $(Q_{n}, \rho_{n})$ with $n$-Auslander-Reiten translation $\tau_{n}$.
The bound Auslander-Reiten quiver of $\mathcal{M}_{n} (D \Gamma^{n} )$ is $(Q_{n+1}, \rho_{n+1}, \tau_{n+1})$.

Iyama proved the following theorem which described the relationship between the bound quivers of $n$-complete algebra $\Gamma^n$ and its cone $\Gamma^{n+1}$ with Auslander-Reiten translations (\cite{i}, Theorem 6.7).

\begin{thm}\label{arq} Let $\Gamma^{n}$ be an $n$-complete algebra which is the cone of an $(n-1)$-complete algebra $\Gamma^{n-1}$ and let $(Q_{n}, \rho_{n})$ be the bound quiver of $\Gamma^{n-1}$ with $n$-Auslander-Reiten translations $\tau_{n} $. Then the bound Auslander-Reiten quiver $(Q_{n+1}, \rho_{n+1})$ with $n$-Auslander-Reiten translation $ \tau_{n+1}$ of $\Gamma^{n}$ is given by the followings:

\begin{itemize}

\item The vertex set $$Q_{n+1,0} =\{ (x, d) | x \in Q_{n, 0}, d\geq 0,  \tau^{d}_{n}x  \neq 0 \},$$ with $$\arr{l}{Q_{n+1,P}=\{(x,d)|x\in Q_{n,0}, \tau_n^{d+1} x = 0 \} \\  Q_{n+1,I}=\{(x,0) | x \in Q_{n,0}\}.}$$

\item The arrow set $Q_{n+1,1}$  consists of two types of arrows:

\begin{itemize}
\item The set of arrows of the first type obtained from  arrows in $Q_{n,1}$: $$\{(\alpha, d): (x, d) \to (y, d)|\alpha: x \rightarrow y \in Q_{n,1}, \mbox{ and  } (x, d), (y, d) \in Q_{n+1,0}\}.$$

\item The set of arrows of the second type obtained by the $n$-Auslander-Reiten translation $\tau_{n}$:
    $$\{ (x, d)_{1}: (x, d) \to (\tau_{n}x, d-1)|(x, d) \in Q_{n+1,0},  d>0\}.$$
\end{itemize}

\item The relation set $$\arr{r}{\rho_{n+1} = \{(r, d)| r \in \rho_{n} \} \cup \{ ( \alpha, d-1)( x, d)_{1} -(\tau_{n}^{-}y, d)_{1}(\tau_n^{-}\alpha, d)\qquad \mbox{\, } \\ ~|~ \alpha: \tau_n x \to y \in Q_{n,1},  d>0 \} ,} $$ here $(r, d)$ is defined as follow: if $p = \alpha_l \cdots \alpha_1$ is a path in $Q_G$, write $(p, t) = (\alpha_l, t) \cdots (\alpha_1, d)$ and if $r =\sum_{v} h_v p_v$ for paths $p_v$ of $Q_G$ and $h_v\in k$, write $(r, d) =\sum_{v} h_v (p_v,d)$.

\item The $n+1$-Auslander-Reiten translation: $$\tau_{n+1}: Q_{n+1,P} \to  Q_{n+1,I},\qquad \tau_{n+1}(x, d)=(x, d+1).$$
\end{itemize}
\end{thm}

\subsection{McKay quiver}
McKay quiver was introduced in 1980 by John McKay for a finite subgroup of the general linear group.
Let $s$ be a positive integer and let $G \subset GL(s,k) = GL(V)$ be a finite subgroup, here $V$ is an $s$-dimensional vector space over $k$.
Let $\{S_i | i =1, 2, \ldots, n\}$ be a complete set of irreducible representations of $G$ over $k$.
For each $S_i$, decompose the tensor product $V \otimes S_i$ as a direct sum of irreducible  representations, write $$V \otimes S_i = \bigoplus_j a_{i,j} S_j, \quad\mbox{  } \quad i = 1, \ldots, n,$$ here $a_{i,j} S_j$ is a direct sum of $a_{i,j}$ copies of $S_j$.
The McKay quiver $Q=Q(G)$ of $G$ is defined the quiver with the vertex set $Q_0 (G)$ the set of (the indices of) the isomorphism classes of irreducible representations of $G$, and there are $a_{i,j}$ arrows from the vertex $i$ to the vertex $j$ in the arrow set $Q_1(G)$.

In \cite{g-m}, Guo and Mart\'inez-Villa proved the McKay quiver of $G$ is the quiver of a pair of skew group algebras $k[V]*G$ and $\wedge V*G$, where $k[V]$ and $\wedge V$ are respectively the symmetric algebra and exterior algebra of the vector space $V$.
Thus there is idempotent $e$ and $\bar{e}$ in $k[V]*G$ and $\wedge V*G$ respectively, such that their basic algebras $V(G) =e k[V]*G e\simeq kQ(G) / (\rho(G)) $ and $\Lambda(G) = \bar{e} \wedge V*G \bar{e} \simeq kQ(G) / (\theta(G)) $ for the path algebra $kQ(G)$ of $Q(G)$ and for some relation sets $\rho(G)$ and $\theta(G)$, respectively.
Since $\wedge V*G $ is a self-injective algebra, we have a {\em Nakayama translation $\nu$} on the McKay quiver $Q(G)$, induced by the Nakayama functor.
We call the bound quiver $(Q_0(G),Q_1(G),\rho(G))$ {\em the bound quiver of $G$.}

It follows from Theorem 2.1 of  \cite{g-m} that the basic algebras $V(G)=V(G)_0+V(G)_1+\cdots$ of $k[V]\ast G$ and the basic algebra $\Lambda(G)= \Lambda(G)_0 +\Lambda(G)_1 +\cdots +\Lambda_{s} $ of $\wedge V\ast G$ are quadratic dual.
They share the same quiver $Q_{G}$, and we have $V(G)_1 = \Lambda(G)_1=kQ_{G,1}$ is the space spanned by the arrows in $Q_{G}$.
The arrows form a basis of $V(G)_1$ and dual basis of $\Lambda(G)_1  = \Hom_k(V(G)_1, k)^{op}$.
That is, we have $\alpha ( \alpha' ) = \left\{ \arr{ll}{ 1 & \alpha = \alpha' \\ 0 & \alpha \neq \alpha'}\right.$, when they are regarded as arrows taken from the quiver of $\Lambda(G)$ and from the quiver $V(G)$, respectively.
We have $\alpha\beta(\alpha'\beta') =\alpha(\alpha') \beta( \beta' )$ in the space $kQ_2$ of paths of length two.
Then relation set $\theta(G)$ and $\rho(G)$ span orthogonal subspaces in $kQ_2$.

In \cite{g-McKay-C}, Guo found that a returning arrow appears from $i$ to $\nu (i)$ at each vertex $i$ in the McKay quiver of $G$ when embedding $GL(V)$ to $SL(V')$ for some $s+1$-dimensional space $V'\supset V$ in some standard manner; and when embedding $G$ into a group $\widetilde{G}$ in $GL(V')$ such that $G = \widetilde{G} \cap SL(V')$ and $\widetilde{G}/ G$ is a cyclic group of order $m$, the McKay quiver of $\widetilde{G}$ is $m$ copies of the McKay quiver of $G$ connected by the arrows induced by the returning arrows. Our main theorem is motivated from this construction of McKay quivers.
We will use trivial extension of self-injective algebras and quadratic duality to give a concrete description of the bound McKay quiver of $\widetilde{G}$.

\subsection{Returning arrows for the trivial extensions}
Recall that {\em the trivial extension $\Lambda\ltimes M$} of an algebra $\Lambda$ by a $\Lambda$-$\Lambda$-bimodule $M$ is the algebra defined on the vector space $\Lambda\oplus M$ with the multiplication defined by $$(a,x)(b,y)=(ab,ay+xb)$$ for $a, b\in \Lambda$ and $x,y \in M$.

Let $M$ be a right $\Lambda$-module and let $\sigma\in$ $\mathrm{Aut}(\Lambda)$.
By $M^{\sigma}$  we denote the right  $\Lambda$-module such that $M^{\sigma}=M$ as $k$-vector space and the right operation of an element $a$ of $\Lambda$ on $M^{\sigma}$ is defined by $ma=m\sigma(a)$, for $m\in M.$
Similarly, $^{\sigma}N$ is defined for a left $\Lambda$-module $N$ and an element $\sigma\in \mathrm{Aut}(\Lambda)$.

Let $\sigma\in \mathrm{Aut}(\Lambda)$.
Denote by $D(\Lambda^{\sigma})$ the dual $\Lambda$-bimodule of $\Lambda^{\sigma}$, and the corresponding trivial extension $T(\Lambda^{\sigma})=\Lambda\ltimes D(\Lambda^{\sigma})$ is called the twisted trivial extension algebra of $\Lambda$.

Fix an integer $l\geq 1$, recall that the bound quiver $Q$ of a graded self-injective algebra $\Lambda$ of Loewy length $l+1$ is a {\em stable translation quiver of Loewy length $l+1$} satisfying the following conditions\cite{g-covering}:

(1). A permutation $\nu$ is defined on the vertex set of $Q$;

(2). The maximal bound paths of $Q$ have the same length $l$;

(3). For each vertex $i$, there is a maximal bound path from $\nu(i)$ to $i$, and there is no bound path of length $l$ from $\nu(i)$ to $j$ for any $j\neq i$;

(4). Any two maximal bound paths starting at the same vertex are linearly dependent.

$\nu$ is called the Nakayama translation of the stable translation quiver $Q$, it is induced by the Nakayama functor.

Let $\Lambda$ be a finite dimensional graded self-injective algebra  given by a stable translation quiver  $(Q,\rho)$ of Loewy length $l+1$ with the Nakayama translation  $\nu$.
By \cite{gyz}, there is a map $\nu$ from $Q_{1}$ to $\bigcup_{i,j}e_{j}\Lambda_{1}e_{i}$ such that $\nu(Q_{1})$ is a basis of $\Lambda_{1}$ which can be extended to a graded automorphism on $\Lambda$.
Actually, $\nu$ is the Nakayama automorphism of $\Lambda$ induced by Nakayama  functor on mod$\Lambda$.
The following quiver $\widetilde{Q} = ( \widetilde{Q}_{0} , \widetilde{Q}_{1})$ is called the returning arrow quiver of $(Q,\rho)$\cite{gyz}, where the vertex set $\widetilde{Q}_{0} = Q_{0}$, the arrow set $\widetilde{Q}_{1} = Q_{1} \cup \{ \alpha_{i}: i \rightarrow \nu (i) ~|~ i\in Q_{0}\}$;

Recall that an automorphism $\sigma$ is called $nice$ if it preserves
idempotents, that is, $\sigma(e)=e$ for all the idempotents of $\Lambda$. The following theorem is proved in  \cite{gyz}.

\begin{thm}\label{gyz} Let $\Lambda$ be a finite dimensional graded self-injective algebra  given by a stable translation quiver  $(Q,\rho)$ of Loewy length $l+1$ with the Nakayama translation $\nu$. If $\sigma$ is a nice graded automorphism of $\Lambda$, then the twisted trivial extension $T(\Lambda^{\sigma})$ is given by the bound quiver $(\widetilde{Q}, \widetilde{\rho}^{\sigma})$, where $\widetilde{\rho}^{\sigma} = \rho \cup \{\alpha_{\nu i } \alpha_{i} ~|~i\in Q_{0}\} \cup \{\alpha_{j} \beta - \sigma^{-1} \nu(\beta)\alpha_{i} ~|~ \beta: i\rightarrow j \in Q_{1}\}$.
\end{thm}

Let $G$ be a finite subgroup of ${\rm Aut}_{k}(\Lambda)$, and let $\Lambda\ast G$ be   the skew group algebra of $\Lambda$, assume that $\sigma\in{\rm Aut}_{k}(\Lambda)$  satisfying $g\sigma=\sigma g$ for $g\in G$,  then $\sigma$  induces a automorphism $\widetilde{\sigma}$ of $\Lambda\ast G$  such that $\widetilde{\sigma}(a\ast g)=\sigma(a)\ast g$ for $a\in \Lambda, g\in G$.
We define $g(a,\varphi)=(g(a),\varphi g^{-1})$ for $(a,\varphi) \in T(\Lambda^{\sigma})$.
Since $g\sigma=\sigma g$, $g$ can be lifted a automorphism of $T(\Lambda^{\sigma})$.
So $G$ is a finite subgroup of ${\rm Aut}_{k} (T( \Lambda^{ \sigma } ))$ and  $T(\Lambda^{\sigma}) \ast G $ makes sense.

\begin{lem}\label{aut}(\cite[Lemma 2.2]{Zheng})
Let $\Lambda$ be a finite dimensional algebra, $G$ a finite subgroup of ${\rm Aut}_{k}(\Lambda)$, assume that $\sigma \in {\rm Aut}_{k}(\Lambda)$ satisfying $g\sigma = \sigma g$ for each $g\in G$. Then we have a  $k$-algebra isomorphism.
$$ T(\Lambda^{\sigma})\ast G \cong T((\Lambda\ast G )^{ \widetilde{\sigma}}).$$
\end{lem}

\begin{rmk}\label{mark} By the proof of Theorem 2.3 in \cite{gyz}, we know that the twisted trivial extension $T ( \Lambda^{ \sigma})$ can also be given by $( \widetilde{Q}, \widehat{ \theta }^{ \sigma})$, where
$$\widehat{\theta}^{\sigma} = \rho \cup \{\alpha_{\tau i} \alpha_{i} ~|~ i\in Q_{0}\} \cup \{ \beta \alpha_{\tau^{-1}(i)} - \alpha_{\tau^{-1}(j)} \nu^{-1} \sigma(\beta) ~|~ \beta: i \rightarrow j\in Q_{1}\}.$$
\end{rmk}

\section{The bound McKay Quiver with Returning  Arrows}
Let $V$ be an $s$-dimensional vector space over $k$, and let $V'$ be an $s+1$-dimensional vector space containing $V$.
Let $G \subset GL(V)$ be a finite subgroup.
$ GL(V) \subset SL(V')$ via the natural embedding $g \to \mat{cc}{g&\\&\hls^{-1} g}$, where $\hls$ is the determinant map.
Let $G'$ be the image of $G$.
The McKay quiver $Q(G')$ is obtained from $Q(G)$ by inserting an arrow $\beta: \nu i \to i$ for each vertex $i$ in $Q(G)$

Let $\Lambda(G)$ and $\Lambda'(G)$  be the basic versions of the skew group algebras $\wedge V*{G}$ and $\wedge V'*{G'}$.
We first study the bound quiver of $\Lambda'(G)$ using trivial extension of graded self-injective algebras.

We need a lemma comparing  certain twisted trivial extensions of a graded algebra with its basic algebra.
Let $A=A_{0}+A_{1}+\cdots +A_{l}$ be a finite dimensional graded algebra such that $A_0$ is semi-simple and $A$ is generated by $A_0$ and $A_1$.
Let $B=B_{0}+B_{1}+ \cdots + B_{l}$  be a basic version of $A$, that is, $B$ is a basic algebra Morita equivalent to $A$.
Then $B_0$ is a direct sum of finite copies of $k$, and generated by $B_0$ and $B_1$.

Let  $\sigma_{0}$ be the graded automorphism of $A$ induced by the map $$\sigma_0: \gamma \mapsto (-1)^{n}\gamma$$  on $A_{1}$.
Then $\sigma_{0}$ induces an automorphism of $B$, which we also denote by $\sigma_{0}$.

\begin{lem}\label{basic}
$T(A^{\sigma_{0}})$ and $T(B^{\sigma_{0}})$ are Morita equivalent.
\end{lem}

\begin{proof}
Take $e\in A_0$ to be an idempotent such that $B=eAe$ is the basic algebra of  $A$, it follows from  Corollary 6.10 of \cite{elements} that
$F={\rm Hom}(Ae,-):{\rm mod}A\rightarrow{\rm mod} B$ is an equivalence  of categories.  Then by Theorem 6.8 of \cite{AK}, there is a projective generator $_{B}P=F(A)=eA$  such that ${\rm End}_{B}(P)\cong A$ as $k$-algebras.
Clearly $P$ is a graded $B$-$A$-bimodule.
By  \cite{FK} there are integers $m$, $n$ and modules $B'$, $P'$ with the following left  $B$-module isomorphisms
$$P^{n}\cong B\oplus B', B^{m}\cong P\oplus P'.$$
Thus we have the left $T(B^{\sigma_{0}})$-module isomorphisms $$(T(B^{\sigma_{0}})\otimes_{B} P)^{n}\cong T(B^{\sigma_{0}})\oplus (T(B^{\sigma_{0}})\otimes_{B} B')$$
and
$$(T(B^{\sigma_{0}}))^{n}\cong (T(B^{\sigma_{0}})\otimes_{B}P)\oplus ((T(B^{\sigma_{0}})\otimes_{B} P')).$$

Thus $T(B^{\sigma_{0}})\otimes_{B}P$ is a projective generator of $\mmod T(B^{\sigma_{0}})$, write it as $\widetilde{P}$.

By Proposition 1.8 and the claim above Proposition 1.13 of \cite{FGR}, we have

$${\rm End}_{T(B^{\sigma_{0}})}(\widetilde{P})\cong {\rm End}_{B}(P)\ltimes {\rm Hom}_{B}(P,D(B^{\sigma_{0}})\otimes P).$$

It follows from  Lemma 5 of \cite{Y}  that $$_{A}D(A)_{A}\cong _{A}{\rm Hom}_{B}(_{B}P_{A},_{B}D(B)\otimes_{B}P_{A})_{A}.$$

Observe  that $$_{A}D(A^{\sigma_{0}})_{A}\cong _{A}{\rm Hom}_{B}(_{B}P^{\sigma_{0}}_{A},_{B}D(B)\otimes_{B}P_{A})_{A}.$$

Define $$\pi:_{A}{\rm Hom}_{B} ( _{B}P_{A}^{\sigma_{0}}, _{B}D(B)_{A} \otimes_{B} P_{A} )_{A} \rightarrow _{A}{\rm Hom}_{B} ( _{B}P_{A}, _{B}D(B^{\sigma_{0}}) \otimes_{B} P_{A})_{A}$$ by $ \pi(f)(p)=\sum_{i}(-1)^{i}f(p_{i}),$ for $f \in _{A}{\rm Hom}_{B}(_{B}P_{A}^{\sigma_{0}}, _{B}D(B)_{A} \otimes_{B} P_{A} )_{A}, p\in P$, $p=\sum_{i}p_{i}$  with $p_{i}\in P_{i}.$

Then $\pi$ is an $A$-$A$-bimodule isomorphism.

Thus we have an algebra isomorphism $T(A^{\sigma_0})\cong {\rm End}_{B}(P)\ltimes {\rm Hom}_{B}(P,D(B^{\sigma_{0}})\otimes P)$. By Theorem 6.7  of \cite{AK},  $T(A^{\sigma_{0}})$ and $T(B^{\sigma_{0}})$ are Morita equivalent.
\end{proof}

Let $\varepsilon$ be a graded automorphism of $\wedge V$ induced by the linear map $$\varepsilon:x \mapsto -x$$ on the vector space $V$.

By Proposition 3.1 of  \cite{AK}, the Nakayama automorphism $\nu$ of $\wedge V$ is induced by the linear map defined by $\nu(x)=(-1)^{n-1}x$ for $x \in V$. Let $\sigma=\varepsilon\nu$, then $\sigma$ is a graded automorphism of $\wedge V$.  For $g=(r_{ij})_{n\times n}\in G$, we have $$ \sigma g(x_{i})=g\sigma(x_{i})=(-1)^{n}\sum^{n}_{j=1}r_{ij}x_{j},$$ for $1\leq i\leq n.$ Therefore $\sigma g=g\sigma$ on $\wedge V$. By the argument above Lemma \ref{aut}, $G$ can be regarded as an automorphism group of $T(\wedge V^{\sigma})$ and $\sigma$ induces an automorphism $\widetilde{\sigma}$ on $\wedge V\ast G$ such  that $$\widetilde{\sigma}(x_{i}\ast g) = \sigma(x_{i}) \ast g = (-1)^{n} (x_{i}\ast g), ~ \widetilde{\sigma}( p \ast g)= \sigma(p) \ast g=p\ast g$$ for $p\in (\wedge V)_{0}=k, 1\leq i\leq n$, and hence $\widetilde{\sigma}$ is a nice automorphism of $\wedge V\ast G$.

We have the following proposition.
\begin{prop}\label{iso}
$\wedge V'\ast G'\cong T((\wedge V\ast G)^{\tilde{\sigma}})$
as algebras.
\end{prop}
\begin{proof}
Note that $D(\wedge V)$ is generated by $(\wedge V)_{n}^{\ast} = D((\wedge V)_{n})$ as $\wedge V$-$\wedge V$-bimodule, and assume that $0\neq \varphi \in (\wedge V)^{\ast}_{n}$.
As an algebra, $T(\wedge V^{\sigma})$ is generated by $(\wedge V)_{0}$ and $(\wedge V)_{1}+(\wedge V)_{n}^{\ast}$, with the additional relation $\varphi^2=0$ and $\varphi x+x \varphi=0$, by Proposition 2.5 of \cite{gyz}.
This shows $T(\wedge V^{\sigma}) \simeq \wedge V'$ and $T(\wedge V^{\sigma})*G \simeq \wedge V'*G'$.
So we get $\wedge V'\ast G'\cong T((\wedge V\ast G)^{\widetilde{\sigma}})$  by Lemma \ref{aut}.
\end{proof}

The following proposition recovers and generalizes Theorem 3.1 of \cite{g-McKay-C}.

\begin{prop}\label{McKay quiver}Suppose that $(Q_{G},\theta_{G})$ is the bound quiver of $\Lambda(G)$, and let $\nu$ be the Nakayama automorphism on $\Lambda(G)$.
The bound quiver $(Q_{G'},\theta_{G'})$ of $\wedge V'\ast G'$ is given by the followings:

1.  $Q_{G'}$ is obtained from $Q_{G}$ by adding an arrow $\beta_{i}: i\rightarrow \nu i$ for each vertex $i\in Q_{G,0}$.

2. $\theta_{G'} = \theta_{G} \cup \{\beta_{\nu i} \beta_{i} ~|~ i \in Q_{G,0}\} \cup \{ \alpha \beta_{\nu^{-1} i} - \beta_{\nu^{-1} j}\nu^{-1} \widetilde{\sigma}(\alpha) ~|~ \alpha: i \rightarrow j \in Q_{G,1}\}.$
\end{prop}

\begin{proof} By Proposition \ref{iso}, $\wedge V'\ast G'\cong T((\wedge V\ast G)^{\widetilde{\sigma}})$.
By Lemma \ref{basic}, we have that $T( \Lambda(G)^{ \widetilde{ \sigma }} )$ is the basic algebra of $T((\wedge V\ast G )^{ \widetilde{ \sigma } })$, and  $ T ( \Lambda( G )^{ \widetilde { \sigma } } ) \cong kQ_{G'}/ ( \theta_{G'})$.
$\Lambda(G)$ is  a finite dimensional graded self-injective algebra with stable translation quiver  $(Q_G,\theta)$ with the Nakayama translation induces by $\nu$.
So by Theorem \ref{gyz},  we get that  $Q_{G'}$ is obtained from $Q_{G}$ by adding an arrow $\beta_{i}: i\rightarrow \nu(j)$ for each vertex $i\in Q_{G,0}$, and by Theorem \ref{gyz} and Remark \ref{mark}, $$ \theta_{G'} = \theta_{G} \cup \{\beta_{ \nu i} \beta_{i} ~|~ i \in Q_{G,0} \} \cup \{\alpha \beta_{\nu^{-1} i} - \beta_{\nu^{-1}j} \nu^{-1} \widetilde{ \sigma } (\alpha) ~|~ \alpha : i \rightarrow j \in Q_{G,1}\}.$$
\end{proof}

Furthermore, we have the following proposition on the bound McKay quiver.
\begin{prop}\label{dualrelations} Let $(Q_{G},\theta_{G})$ be the bound quiver of $\wedge V\ast G$, then the relation set of the bound McKay quiver $(Q_{G'}, \rho_{G'}) $ of $G'$ is given by
$$\rho_{G'}  = \rho_{G} \cup \{\alpha\beta_{\nu^{-1}(i)} - \beta_{\nu^{-1}(j)} \nu^{-1} (\alpha) ~|~ \alpha: i \rightarrow j \in Q_{G,1} \}. $$
\end{prop}
\begin{proof}
Let $Q_{G,2}$ and $Q_{G',2}$ be the sets of the paths of length $2$ of $Q_{G}$ and of $Q_{G'}$. respectively.
Since $Q_{G,1}$ is a subset of $Q_{G',1}$, we have $\rho(G) \subset \rho(G')$, $\theta(G)\subset \theta(G')$ and $kQ_{G,2} \subset kQ_{G',2}$.
Since $V(G)$ and $\Lambda(G)$ are quadratic dual, and $V(G)$ and $\Lambda(G)$ are quadratic dual, $\rho(G)$ and $\theta(G)$ span orthogonal subspaces in $kQ_{G,2}$, and $\rho(G')$ and $\theta(G')$ span orthogonal subspaces in $kQ_{G',2}$.
For each pair $i, j$ of vertices, the subsets $e_j \theta_{ij}( G) e_i$ and $e_j \rho(G) e_i$  span of orthogonal subspaces in $e_j kQ_{G,2} e_i$, and the subsets $e_j \theta_{ij}( G') e_i$ and $e_j \rho(G') e_i$  span of orthogonal subspaces in $e_j kQ_{G',2} e_i$.

Fix a vertex $i$, the set the paths of length $2$ from $i$ to $j$  of $Q_{G'}$ is $$e_jQ_{G',2}e_i =e_jQ_{G,2}e_i\cup \{\beta_{\nu i}\beta_i\} \cup \{\beta_{\nu^{-1}j} \alpha, \nu(\alpha)' \beta_i | \alpha \in e_{\nu^{-1}j} Q_1 e_i \} $$ if $j= \nu^2(i)  $, and is $$e_jQ_{G',2}e_i =e_jQ_{G,2}e_i\cup \{\beta_{\nu (i)}\beta_i\} \cup \{\beta_{\nu^{-1}(j)} \alpha, \nu(\alpha)' \beta_i | \alpha \in e_{\nu^{-1}(j)} Q_1 e_i \} $$
if $j\neq \nu^2(i) $.
Thus $e_j \rho(G')e_i$ span an orthogonal subspace of $e_j\theta(G') e_i$ since $$\arr{lll}{&& (\alpha\beta_{\nu^{-1}(i)}- \beta_{\nu^{-1}(j)} \nu^{-1}\widetilde{ \sigma } (\alpha)) (\alpha\beta_{\nu^{-1}(i)}+ \beta_{\nu^{-1} (j)}\nu^{-1} \widetilde{ \sigma } (\alpha))\\ &=& \alpha(\alpha)\beta_{\nu^{-1}(i)}(\beta_{\nu^{-1}(i)})- \beta_{\nu^{-1}(j)}(\beta_{\nu^{-1}(j)}) \nu^{-1}\widetilde{ \sigma } (\alpha)(\nu^{-1}\widetilde{ \sigma } (\alpha))\\&=&0}.$$
Since $\widetilde{ \sigma } (\alpha) =-\alpha$ for all $\alpha \in Q_{G,1}$, this shows that $$\rho_{G'}  = \rho_{G} \cup \{\alpha\beta_{\nu^{-1}(i)} - \beta_{\nu^{-1}j} \nu^{-1} \alpha ~|~ \alpha: i \rightarrow j \in Q_{G,1} \} $$ is the relation set for $V(G')$.
\end{proof}

\section{Bound McKay Quivers for Cyclic Extension}

By embedding $GL(V)$ into $SL(V')$, we obtain returning arrows in the McKay quiver of a group.
We need smash product construction for the algebra to change the returning arrows into connecting arrows between the copies of the original McKay quiver.
These connecting arrows are the arrows of type $(x,d)_1: (x,d) \to (\tau_n x, d-1)$ in the Auslander-Reiten  quiver of the cone (Theorem \ref{arq}).
Such construction is realized by a direct product of group $G'$ with a cyclic group.
Let $H$ be a finite group,  and let  $\Lambda$ be an $H$-graded algebra.
Recall that the smash product $\Lambda \#kH^*$ is the free $\Lambda$-module with basis $H^{\ast}=\{\delta_{g} ~|~ g\in H\}$,  and the multiplication is defined by $$a\delta_{g} \cdot b \delta_{h} =a b_{gh^{-1}} \delta_{h},$$
for $a, b\in \Lambda$, where $b_{gh^{-1}}$ is the homogeneous component of degree $gh^{-1}$ of $b$ in $\Lambda$.

Let $\xi_{m}\in k$ be a primitive $m$-th root of unity and let $$\omega_{m}=\left(   \begin{array}{cccc}    1 &  &  &  \\     & \ddots &  &  \\     &  & 1 &  \\     &  &  & \xi_{m} \\  \end{array} \right)_{(s+1)\times(s+1)}$$
Let $C_m=(\xi_{m})$ be the cyclic subgroup of $GL(V')$ generated by $\omega_{m}$ and let $\widetilde{G}=GC_m = G\times C_m$, then $\widetilde{G}\cap SL(V')=G$.
Let $\widehat{C_{m}}={\rm  Hom}(C_{m},k^{\ast})\cong \mathbb{Z}_{m}$ be the dual group of $C_{m}$.
For $\lambda_{i}\in {\rm  Hom}(C_{m},k^{\ast})$ such that $\lambda_{i} ( \omega_{m} ) = \xi_{m}^{i}$, then $\lambda_{i} =\lambda_{ 1 }^{ i } $ and $\widehat{C_{m}}$ is the cyclic group generated by $\lambda_{1}$.

Let  $V'=V+V_{1}$ is a direct sum of subspaces,  and $V_{1}$ is an one dimensional subspace of spanned by $y$.
The action of $C_m$ on $V'$ induce a $C_m$-grading on $\wedge V'\ast G'$ with such that $a\in (\wedge V'\ast G')_{\omega^i}$ if and only if $\omega(a) = \xi^i a$.
Thus $\wedge V'\ast G' = (\wedge V'\ast G')_{\omega_0}+(\wedge V'\ast G')_{\omega_1}$, with $(\wedge V'\ast G')_{\omega_0}= \wedge V \ast G'$ and $(\wedge V'\ast G')_{\omega_1}= y(\wedge V \ast G')$,
is a $C_m$-graded algebra.
Now it follows from \cite{p} that
$$\wedge V'\ast \widetilde{G}\cong \wedge V'\ast (G'\times C_{m})\cong (\wedge V'\ast G')\ast C_{m}, $$
and by \cite{ms}, $(\wedge V'\ast G')\ast C_{m}\cong (\wedge V'\ast G')\#k\widehat{C_{m}}^{\ast}.$

Clearly, if $e_i$ is a primitive idempotent of $\wedge V'\ast G'$, then for each $j\in \mathbb Z/mZ$,  $e_i\delta_j$ is a primitive idempotent of $(\wedge V'\ast G')\#k\widehat{C_{m}}^{\ast}$, and we have the following lemma.
\begin{lem}\label{new} Suppose that $e_i$ and $e_i'$ are primitive  idempotents of $\wedge V'\ast G'$, and $ e_i(\wedge V'\ast G')\cong e_i'(\wedge V'\ast G')$ as right $\wedge V'\ast G'$-modules, then for $j\in\mathbb{Z} / {m} \mathbb{Z}$, $e_i\delta_{j}((\wedge V'\ast G')\#k\widehat{C_{m}}^{\ast})\cong e_i'\delta_{j}((\wedge V'\ast G')\#k\widehat{C_{m}}^{\ast})$ as right $\wedge V'\ast G'\#k\widehat{C_{m}}^{\ast}$-modules.
\end{lem}
\begin{proof}
If $\varphi: e_i(\wedge V'\ast G')\to e_i'(\wedge V'\ast G')$ is an isomorphism as right $\wedge V'\ast G'$-modules with $\varphi(e_i)=e_i'b$ for $b\in\wedge V'\ast G'$,  then $$\varphi(e_i)_{0}\delta_{t}=e'_ib_{0}\delta_{t}=e_i' \delta_{t} b \delta_{t} \in e'\delta_{t}((\Lambda V'\ast G')\#k\widehat{C_{m}}^{\ast}).$$

Set $\varphi' ( e_i \delta_{t} ) = \varphi( e_i )_{0} \delta_{t}$, we get right $\Lambda V'\ast G' \# k\widehat{ C_{m} }^{\ast}$-module isomorphism $$ \varphi': e_i \delta_{t} ((\Lambda V'\ast G') \# k \widehat{ C_{m} }^{ \ast}) \rightarrow e_i'\delta_{t}( ( \Lambda V'\ast G') \# k \widehat{ C_{m} }^{\ast}).$$
\end{proof}

Let $e = e_{1}+e_{2}+\ldots +e_{r}$ be a sum of orthogonal primitive idempotents in $\Lambda V'\ast G'$ such that $e \Lambda V' \ast G'e = \Lambda(G')$ is its basic algebra.
Let $\widetilde{e}= \sum^{m}_{j=1} \sum^{r}_{i=1} e_i \delta_{j} =\sum^{m}_{j=1} e\delta_{j}$, then $\tilde{e}$ is an idempotent in $(\Lambda V'\ast G')\#\widehat{C_{m}}^{\ast}$.

We also have the following lemma.
\begin{lem}\label{equal}
$$\Lambda(G') \# k \widehat{C_{m}}^{\ast} = \widetilde{e} ((\wedge V'\ast G')\#k\widehat{C_{m}}^{\ast})\widetilde{e}.$$
\end{lem}
\begin{proof} Clearly, $\Lambda(G') \# k \widehat{C_{m}}^{\ast} \subseteq \widetilde{e} ((\wedge V'\ast G' ) \# k \widehat{ C_{m} }^{\ast}) \widetilde{e}.$
For any $\sum^{m}_{t=1}a_{t}\delta_{t} \in (\wedge V' \ast G' ) \# k \widehat{C_{m}}^{\ast},$ with $a_{t}\in\wedge V'\ast G'$,   write $[a_{t}]_j$ for the $\omega_j$ component of $a_t$, we have that
$\widetilde{e}(\sum^{m}_{t=1}a_{t}\delta_{t}) \widetilde{e}
= \sum^{m}_{t=1}\sum^{m}_{j=1} \sum^{m}_{j'=1} e \delta_{j} \cdot (a_{t}\delta_{t} \cdot e \delta_{j'}) = \sum^{m}_{t=1} e \sum^{m}_{j=1} [a_{t}]_{j-t} e \delta_{t} = \sum^{m}_{t=1} ( e a_{t} e ) \delta_{t}.$
Thus $$\Lambda(G') \# k \widehat{C_{m}}^{\ast} \supseteq \widetilde{e} ( ( \wedge V' \ast G') \# k \widehat{C_{m}}^{\ast} ) \widetilde{e},$$ and hence
$$ \widetilde{e}  ( ( \wedge V' \ast G') \# k \widehat{C_{m}}^{\ast})\widetilde{e}  = \Lambda(G') \# k \widehat{C_{m}}^{\ast}.$$
\end{proof}

Now we describe the bound quiver of $\Lambda(G') \# k \widehat{ C_{m} }^{\ast}$ using the bound quiver of $\Lambda(G)$.

\begin{thm} \label{a} Let $(Q_{G}, \theta_{G})$ be the bound  quiver of $\Lambda(G)$ with the Nakayama automorphism $\nu$, then the bound quiver $(Q_{\widetilde{G}}, \theta_{ \widetilde{G} } )$ of $\Lambda( \widetilde{G})$ is the quiver defined by the following data:

1. The vertex set $Q_{\widetilde{G},0} = Q_{G,0} \times \mathbb{Z}/{m}\mathbb{Z}$.

2. The arrow set $Q_{ \widetilde{G}, 1} = \bigcup_{t \in \mathbb{Z}/{m}\mathbb{Z}} (\{ (\alpha,t) : (i,l) \rightarrow (j,t) ~|~ \alpha: i \rightarrow j\in Q_{G,1} \} \cup \{ (\beta_{i} , t):(i,t) \rightarrow (\nu i, t + 1) ~|~ i \in Q_{G,0}, \})$.

3. The relation set  $$\arr{cl}{\theta_{ \widetilde{G} } &= \bigcup_{t \in \mathbb{Z}/{m}\mathbb{Z}} ( \{(\varrho, t) ~|~ \varrho \in \theta_{G}\} \cup \{ ( \beta_{\nu(i)}, t+1) (\beta_{i},t) ~|~i \in Q_{G,0}\}\\ &\qquad \cup \{ (\alpha, t) (\beta_{\nu^{-1}i}, t-1) - (\beta_{\nu^{-1} j}, t-1) (\nu^{-1}\tilde{\sigma} \alpha, t-1)  ~|~ \alpha: i \rightarrow j \in Q_{G, 1}\}),}$$ here $(\varrho, t)$ is defined as follows: if $p=\alpha_l\cdots\alpha_1$ is a path in $Q_G$, write $(p,t)=(\alpha_l,t)\cdots(\alpha_1,t)$ and if $\varrho =\sum_{v} h_v p_v$ for paths $p_v$ of $Q_G$ and $h_v\in k$, write $(\varrho,t) =\sum_{v} h_v (p_v,t)$ for all $t\in \mathbb Z/m\mathbb Z$.

4. The Nakayama translation of $\Lambda(\widetilde{G})$  is defined by $\nu(i, t)=(i, t-1),$ for $(i,t)\in Q_{ \widetilde{G}, 0} .$
\end{thm}
\begin{proof}
By Lemma \ref{equal}, $(Q_{\widetilde{G}}, \theta_{ \widetilde{G} } )$  is the bound quiver $\Lambda(G')\# \widehat{C_{m}}^{\ast}$.

Since $\Lambda(G')$ is a graded self-injective algebra of Loewy length $s+2$, it is an $s$-translation algebra and $(Q_{G'}, \theta(G')$ is an admissible $s$-translation quiver.
The theorem follows from  Proposition \ref{McKay quiver} and Proposition 5.5 of \cite{g-tran}.
\end{proof}

We have the following corollary.

\begin{coro}\label{McKay quiver 2}
Let $(Q_{G},\rho_{G})$ be the bound McKay quiver of $G$ with the Nakayama automorphism $\nu$, then the bound quiver  of $ \widetilde{G}$ is $(Q_{\widetilde{G}}, \rho_{\widetilde{G}})$, where $Q_{\widetilde{G}}$ is as in the above theorem and $$\arr{r}{\rho_{\widetilde{G}}=\bigcup_{t\in \mathbb Z/m\mathbb Z}(\{(\varrho, t) ~|~ \varrho \in \rho_{G}\} \cup \{ (\alpha, t) (\beta_{\nu^{-1} i},t-1) - ( \beta_{\nu^{-1} j}, t-1 ) (\nu^{-1} \alpha, t-1) \\ ~|~ \alpha: i \rightarrow j \in Q_{G, 1},l\in \mathbb{Z}_{m}\}). }$$
\end{coro}
\begin{proof}
By Theorem  \ref{a}, we get the bound quiver $(Q_{\widetilde{G}}, \theta_{\widetilde{G}})$ of  $\Lambda( \widetilde{G})$.
Note that $V(\widetilde{G})$ and $\Lambda(\widetilde{G})$ are quadratic dual, similar to the proof of Corollary \ref{dualrelations}, one sees that $\rho_{\widetilde{G}}$ spans the orthogonal subspace of $\theta_{\widetilde{G}}$ in $kQ_{\widetilde{G}, 2}$.
Thus  $(Q_{\widetilde{G}}, \rho_{\widetilde{G}})$ is the bound McKay quiver of $\widetilde{G}$.
\end{proof}

\section{Main theorems}
In this section we shall prove our main theorem.

Let $(Q,\rho)$ be a bound quiver and let $Q'$ be a full bound  subquiver of $Q$ in the sense that the arrows from $i$ to $j$ in $Q'$ are the same as in $Q$ provided that both $i,j$ are vertices in $Q'$.
For each element $a = \sum_{v \in J_a} h_vp_v \in kQ$, where $p_v$ are paths in $Q$ and $h_v\in k$.
Let $J'_a=\{v\in J_a| \mbox{the vertices of } p_v  \mbox{ are all in }Q'\}$ and write $a_{Q'}=\sum_{v \in J'_a} h_vp_v $.
We say that $a_{Q'}$ is the component of $a$ in $Q'$.
For each $v\in J_a\setminus J'_a$, we have that $p_v$ passes through a vertex outside $Q'$, thus we have  $p_v=q_ve_{i_v}q_v'$ for $i_v$ in $Q_0\setminus Q'_0$ in $kQ$.
Let $\rho'= \{a_{Q'}| a\in \rho\}$ be the relation set on $Q'$ induced by $\rho$ and call  $(Q',\rho')$ {\em a full bound  subquiver of  $(Q, \rho)$}.

Let  $Q'$  and $Q''$ be finite quivers, write $s, t$ for the maps from the arrow set to the vertex set sending an arrow to its source and its target, respectively.
So we have  $\alpha: s(\alpha) \to t(\alpha)$ for an arrow $\alpha$.
A pair of maps  $\omega=(\omega_{0},\omega_{1})$, where $\omega_{0}: Q'_{0}\rightarrow Q''_{0}, \omega_{1}: Q'_{1} \rightarrow Q''_{1}$, is called  a {\em quiver embedding } if it satisfies the following conditions:

(i). $s\omega_{1}(\alpha)=\omega_{0} s(\alpha),  t\omega_{1}(\alpha)=\omega_{0} t(\alpha)$ for $\alpha\in Q_{1}'$.

(ii). $\omega_{0}$ and $\omega_{1}$ are injections.

In this case, we  write $\omega(Q')\subseteq Q''$.

Obviously, $\omega$ induces a algebra monomorphism $\omega : kQ' \rightarrow kQ$.

Let  $(Q', \rho')$ , $(Q, \rho)$ be bound quivers, and $\omega$ an algebra monomorphism $\omega : kQ' \rightarrow kQ$, write $\omega(\rho')= \{\omega(x)| x \in \rho' \}$.
We introduce the following definition.

\begin{defn}
The bound quiver $(Q', \rho')$ is said to be a {\em truncation} from the bound quiver $(Q,\rho)$, if there is a
quiver embedding  $\omega: Q' \rightarrow Q$ satisfying that $\omega(\rho')$ is induced by $\rho$.
\end{defn}

In this case, the quiver embedding $\omega$ is called a bound quiver truncation of $(Q', \rho') $ from $(Q, \rho)$, and the algebra $\Lambda'=kQ'/(\rho')$ is called a { truncation} of the algebra $\Lambda=kQ/(\rho)$.

We have the following proposition.

\begin{prop} Let $\Lambda'$ and $ \Lambda$ be algebras over $k$ given by the bound quivers $(Q', \rho')$ and $(Q, \rho)$, respectively.
If $\Lambda'$ is a truncation of $\Lambda$, there is a subset $E$ of the orthogona primitive idempotents of $\Lambda$ such that $\Lambda' \cong \Lambda/(E)$, where $(E)$ is the ideal generated by $E$.
\end{prop}
\begin{proof}
Since $\Lambda'$ is a truncation of $\Lambda$, $(Q', \rho')$ is  a truncation from $(Q, \rho)$, induced by some quiver embedding  $\omega$ and $\Lambda' = kQ'/(\rho')$.
Let $\widetilde{E}$ be the preimage of $E$ in $kQ$, that is, the set of trivial paths in $kQ$ which induce primitive idempotents in $\Lambda$.
Then we have an isomorphism $kQ/(\widetilde{E}) \simeq kQ'$, write $\phi: kQ \to kQ'$ for the homomorphism which induces the isomorphism, then $\Ker \phi =(\widetilde{E})$.
Let $\pi': kQ' \to \Lambda'$ be the canonical homomorphism, then $\pi' \phi: kQ \to \Lambda'$ is an epimorphism and we have that $\Ker \pi' \phi$ is generated by  $E\cup \omega(\rho')$, that is.
$\Ker \pi' \phi = (E, \omega(\rho')) = (E, \{a_{\omega(Q') } | a \in \rho\}) = (E, \rho)$, since $qe_iq'$ in $(E)$ if $i\in E$.
In particular, $(\rho)\subset \Ker \pi \phi $ and the map $\pi' \phi = \psi \pi $ factor through the canonical homomorphism $\pi: kQ \to kQ/(\rho) = \Lambda$.
Thus $\psi : \Lambda\to \Lambda'$ is an epimorphism and $\Ker \psi = \pi(\widetilde{E}) =(E)$.
That is $\Lambda' \cong \Lambda/(E)$.
\end{proof}

Let $\Gamma\sk{n}$ be an $n$-complete algebra over $k$ which is the cone of an  $n-1$-complete algebra and let $(Q_n, \rho_n)$ be its bound quiver with $n$-Auslander-Reiten translation $\tau_n$.
Let $G\subset GL(n, k)$ be a finite subgroup and let $(Q_G, \rho_G)$ the McKay quiver of $G$ with the Nakayama translation $\nu_{G}$.
We say that the bound quiver of $\Gamma\sk{n}$ is {\em a truncation of the bound McKay quiver of $G$} if there is a truncation $\omega\sk{ n }$ of $(Q_n, \rho_n)$ from $(Q_G, \rho_G)$ such that $$\omega\su{n}{0}\tau_n = \nu_G \omega \su{ n }{ 0 }.$$

Now we state and prove our main theorem.

\begin{thm}\label{man}
Let $\Gamma\sk{n}$ be an $n$-complete algebra over $k$ which is the cone of an  $n-1$-complete algebra and let $\Gamma\sk{n+1}$ be the cone of $\Gamma\sk{n}$. Assume that there is a finite subgroup $G\subset GL(n, k)$ such that the bound quiver of $\Gamma\sk{n}$ is a truncation of the McKay quiver of $G$.

Then there is $m$, and a finite subgroup $\widetilde{G} $ of $GL(n+1, k)$ such that the bound quiver of $\Gamma\sk{n+1}$ is a truncation of the McKay quiver of $\widetilde{G}$.
\end{thm}
\begin{proof}
Let $(Q_n, \rho_n)$ and $(Q_{n+1}, \rho_{n+1})$ be the bound quiver of $\Gamma\sk{n}$ and $\Gamma\sk{n+1}$, respectively.
Then $(Q_{n+1}, \rho_{n+1})$ is constructed from $(Q_n, \rho_n)$ by Iyama as in Theorem \ref{arq}.

Let $m $ be a sufficient large integer such that $\tau_n^m i=0$ for all vertex $x\in Q_{n,0}$.
Let $G'$ be the subgroup of $SL(n+1,k)$ which is the image of $G$ under the canonical embedding $g \to \mat{cc}{g \\ & \hls^{-1} (g)}$, and let $C_m$ be the subgroup of $GL(n+1,k)$ generated by $\mat{cccc}{1\\&\ddots \\ &&1 \\ &&&\zeta_m}$, where $\zeta_m$ is a primitive $m$th root of the unity.
Let  $\widetilde{G} =G'\times C_m$ be a finite subgroup of $GL(n+1, k)$.
The bound McKay quiver of $\widetilde{G}$ is given in Corollary \ref{McKay quiver 2}.

Let $\omega$ be a truncation of the bound quiver of $\Gamma\sk{ n }$ from the McKay quiver of $G$.
Now we define a quiver embedding $\wto =(\wto_0, \wto_1)$ from $(Q_{n+1}, \rho_{n+1})$ to the bound McKay quiver $( Q_{ \widetilde{ G } }, \rho_{ \widetilde{ G }})$:
Define $\wto_0$ on the vertex set by $$\wto_0(x,d) = ( \omega_0 ( x ), -d)$$ for $(x, d)\in Q_{n+1, 0}$.
For an arrow $\alpha: x \to y$ in $Q_{n,1}$ and $d\ge 0$ with $\tau_n^d x \neq 0 \neq \tau_n^d y$, $(\alpha, d)$ is an arrow of the first type in $Q_{n+1, 1}$, we define $\wto_1 (\alpha, d) = (\omega_1(\alpha), -d)$.
For each vertex $x$ and $d\ge 0$ with $\tau_n^{d+1} \neq 0$, $(x, d)_1$ is an arrow of the second type in $Q_{n+1,1}$, we define $\wto(x, d)_1 = (\beta_{\omega_0(x)}, -d)$.
One sees easily that $\wto$ is a quiver embedding from $Q_{n+1}$ to $Q_{\widetilde{G}}$.

For a relation of the form $(r,d)\in \rho_{n+1}$ with $r\in \rho_{n}$,  $\omega(r)  $ is the component in $\omega(Q_n)$ of some relation $\varrho \in \rho_{G}$.
Thus $\wto(r,d) = (\omega(r), -d)$ is the component in $\omega ( Q_{n+1} )$ of the relation $(\varrho,d) \in \rho_{\wt{G}}$.
For a relation of the form  $( \alpha, d-1)( x, d)_{1} -(\tau_{n}^{-}y, d)_{1}(\tau_n^{-}\alpha, d)$ with $\alpha: \tau_n x \to y \in Q_{n,1}$ and $d>0$, write $\omega(\alpha) =\alpha$, $-d=t-1$ then $$\arr{ccl}{&& \wto(( \alpha, d-1)( x, d)_{1} -(\tau_{n}^{-}y, d)_{1} (\tau_n^{-} \alpha, d))\\ & =& (\omega( \alpha), -d+1)\wto( x, -d)_{1} -\wto (\tau_n^{-}y, -d)_{1}(\omega\tau_n^{-}(\alpha), -d)\\ &=&( \alpha, t) ( \beta_x, t-1) -(\beta_{\nu^{-1}y}, t-1) ( \nu^{-1} \alpha, t-1),}$$
which is in $\rho_{\wt{G}}$.
By comparing Theorem \ref{arq} with Theorem \ref{a} and Corollary \ref{McKay quiver 2}, we see that $\wto$ is a truncation of the bound quiver of $\Gamma^{n+1}$ from the McKay quiver of $\wt{G}$.
\end{proof}

Fix an integer $s \geq 1$, Iyama describes a family of $(n-1)$-Auslander absolute $n$-complete algebras $T^{(n)}_{s}(k)$ using quivers with relations, for $n=1,\cdots s$ \cite{i}.
Let $A_{s}$ be a quiver  with $A_{s,0}=\{1,\ldots,s\},$ $A_{s,1}=\{a_{i}:i\rightarrow i+1~|~i=1,\ldots,s-1\}$, then $A_{s}$ is a linear oriented Dynkin quiver of type $A$.
Suppose  that  $T^n_{s}(k)\cong (kQ^{n})/(\rho_{n})$. As an immediate application of    Theorem \ref{man},  we give a new proof of the first part of Theorem 4.6 of \cite{g-McKay}.

\begin{example}
Let $T^{n}_{s}(k)$ be an $(n-1)$-Auslander absolute $n$-complete algebra  for
$n\geq 1$,  then there exists a finite abelian group $G_{n}\subset GL(n,k)$, such that the  bound quiver  of  $T^{n}_{s}(k)$ is a
truncation of the bound McKay quiver $(Q_{G_{n}}, \rho_{G_{n}})$ of $G_{n}$.
\end{example}
\begin{proof}  For $n=1$, $T_{s}^{1}(k)\cong kA_{s}.$ Let $G_{1}=\mathbb{Z}_{r_{1}}\subseteq GL(1,k)$, where $r_{1}$ is sufficiently large. By  Proposition 3.2 of \cite{g-McKay}, $Q_{G_{1},0}=\mathbb{Z}_{r_{1}}$ and $Q_{G_{1},1}=\{\alpha_{i}:i\rightarrow i+1~|~i\in \mathbb{Z}_{r_{1}}\}$. We define the map: $\varsigma_{1}:A_{s}\rightarrow Q_{G_{1}}$ as follows: $\varsigma_{1}(i)=i$, $\varsigma_{1}(a_{i})=\alpha_{i}$, for $i\in A_{s,0},a_{i}\in A_{s,1}$. It's straightforward to check that $\varsigma_{1}$  is a quiver homomorphism, and observe that $\rho_{A_{s}}=\rho_{G_{1}}=\emptyset$. Then $A_{s}$ is a truncation of $(Q_{G_{1}},\rho_{G_{1}})$ induced by   $\varsigma_{1}$, and $\nu_{1}\varsigma_{1}(i)=i-1=\varsigma_{1}\tau_{1}(i)$, $\nu_{1}\varsigma_{1}(a_{i})=\alpha_{i-1}=\varsigma_{1}\tau_{1}(a_{i})$, for $i\in A_{s,0},a_{i}\in A_{s,1}$, where  $\nu_{1}$ is the  Nakayama translation defined on the McKay quiver $Q_{G_{1}}$.

By Theorem 1.19 of \cite{i}, $T^{n}_{s}(k)$ is the cone of $T^{n-1}_{s}(k)$, and $T^{n-1}_{s}(k)$ is an $(n-1)$-complete algebra. Then the claim follows  by induction on $n$  and Theorem \ref{man}.
\end{proof}

\end{document}